\documentclass[11pt]{article}
\usepackage{color,colortbl} 
\usepackage[utf8]{inputenc}
\usepackage{amssymb,amsmath,amsfonts,amsthm,moresize,xfrac,nicefrac,accents,enumerate}
\usepackage[all,cmtip]{xy}
\usepackage{newtxmath}
\usepackage{mathrsfs}
\usepackage{url}
\usepackage{xspace}
\usepackage{tabularx}
\usepackage{booktabs}
\usepackage{tikz-cd}
\usepackage{sectsty}
\usepackage{listings}
\usepackage{tcolorbox}
\tcbuselibrary{breakable}

\definecolor{dkgreen}{rgb}{0,0.6,0}
\definecolor{gray}{rgb}{0.5,0.5,0.5}
\definecolor{mauve}{rgb}{0.58,0,0.82}
	
\definecolor{ltred}{rgb}{1,.88,.88}
\definecolor{ltyellow}{rgb}{.88,.88,1}
\definecolor{ltgreen}{rgb}{.88,1,.88}
\definecolor{ltltgreen}{rgb}{.94,1,.94}
\definecolor{ltpink}{rgb}{1,1,.95}
\definecolor{ltpurple}{rgb}{1,.95,1}
\definecolor{bkgd2}{rgb}{1,1,.92}
\definecolor{v1color}{rgb}{1,.8,.86}
\definecolor{v2color}{rgb}{.8,.88,.92}
\definecolor{v3color}{rgb}{.8,1,.8}

\allsectionsfont{\normalfont\rmfamily\bfseries}

\usepackage{geometry,parskip,cmap,dcolumn}
\allowdisplaybreaks
\usepackage{multicol}
\newtheorem{theorem}{Theorem}
\newtheorem{corollary}{Corollary}

\newcommand{\poincare}{Poincar\'{e}\xspace}
\newcommand{\eye}{\mathbf{I}}

\newcommand{\db}[1]{\mb{e}^*_{#1}}
\renewcommand{\emph}[1]{\textit{#1}}

\newcommand{\mbf}[1]{\mathbf{#1}}
\newcommand{\quot}[1]{``#1''}
\newcommand{\e}[1]{\mb{e}_{#1}}
\newcommand{\myboldhead}[1]{\vspace{0in}\hspace{-0.21in}\textbf{#1.}}

\newcommand\R[1]{\ensuremath{\mathbb{R}^{#1}\xspace}}
\newcommand\RP[1]{\ensuremath{\mathbb{R}{P^{#1}}\xspace}}
\newcommand{\ej}[1]{\mbf{e}^{#1}}
\newcommand{\pdclal}[3]{P({\mathbb{R}^*_{#1,#2,#3}}\xspace)}
\newcommand{\pclal}[3]{P({\mathbb{R}_{#1,#2,#3}}\xspace)}

\usepackage{sectsty}
\def \mb {\mathbf}
\newcommand{\MP}{P\xspace}
\newcommand{\CR}{H\xspace}
\newcommand{\PN}{J\xspace}
\newcommand{\RL}{R\xspace}

\begin{document}

\title{A bit better: Variants of duality in geometric algebras with degenerate metrics}
\author{Charles Gunn, Ph. D.\\ Raum+Gegenraum\\ Falkensee, Germany}
\maketitle

\begin{abstract}

Multiplication by the pseudoscalar $\mathbf{I}$ has been traditionally used in geometric algebra to perform non-metric operations such as calculating coordinates and the regressive product.  In algebras with degenerate metrics, such as euclidean PGA $P(\mathbb{R}^*_{3,0,1})$, this approach breaks down, leading to a  search for non-metric forms of duality. The article compares the dual coordinate map $J: G \rightarrow G^*$, a \emph{double algebra} duality, and Hodge duality $H: G \rightarrow G $, a \emph{single algebra} duality for this purpose. While the two maps are computationally identical, only $J$ is coordinate-free and provides direct support for \emph{geometric duality}, whereby every geometric primitive appears twice, once as a point-based and once as a plane-based form, an essential feature not only of projective geometry but also of euclidean kinematics and dynamics. The article concludes with a proposed duality-neutral software implementation, requiring a single bit field per multi-vector. 

\end{abstract}


\section{History}
\label{sec:his}
Until recently, almost all literature on geometric algebras assumed that the metric was non-degenerate.  A typical justification for this assumption can be found in \cite{hlr01}, Sec. 13.8:
\begin{quote}
Any degenerate algebra can be embedded in a non-degenerate algebra of larger dimension, and it is almost always a good idea to do so.  Otherwise, there will be subspaces without a complete basis of dual vectors, which will complicate algebraic manipulations.
\end{quote}
In this context, the \emph{dual} of a multi-vector $X$ in a geometric algebra $G$ is meant the orthogonal complement, obtained by pseudoscalar multiplication $P: P(X) := XI$ (or some closely-related variant).  When $I^2\neq0$, $P$ is a grade-reversing map that is a vector space isomorphism on each $\bigwedge^k$. $P$ can then be used to carry out crucial computations such as the calculating coordinates and the regressive product (see Secs. \ref{sec:cc} and \ref{sec:regpro} below).  When the metric is degenerate ($I^2=0$), however, $P$ is no longer an isomorphism and these algorithms cease to function.

At about the same time as the above citation, Jon Selig (\cite{selig00}) announced the existence of a geometric algebra for $n$-dimensional euclidean space $\mathbf{E}^{n}$ with degenerate signature $(n,0,1)$, that is, one basis vector that squares to 0. This algebra was worked out in more detail by the current author in \cite{gunnthesis},\cite{gunn2011},\cite{gunn2017a}, and \cite{gunn2017b}. Selig and Gunn suggested different approaches to avoid the use of $P$, based on the common observation that the calculations in question are intrinsically non-metric, and that the use of $P$ is a matter of convenience, not necessity.  


In the meantime a further alternative to $P$ has been proposed. Combined with the original $P$, these four variants present a daunting complexity to the newcomer to the field. This article gives an overview of these different approaches, point out differences and similarities, and propose a simple software design allowing all to co-exist peacefully to a large extent. The article also proposes a terminology to allow the different forms of duality to be conceptually differentiated.


\section{Outline}

The article begins with reviewing some prerequisites for understanding its contents in Sec. \ref{sec:prel}.
It then turns to consider metric and non-metric strategies for calculating coordinates in Sec. \ref{sec:cc}. This discussion leads to the introduction of the \textit{right complement} map in Sec. \ref{sec:rcm}, which plays a large role in what follows.
Sec. \ref{sec:regpro} addresses strategies to calculate the regressive product. This features four different approaches: the metric one (valid only for non-degenerate metrics), the shuffle product (from the theory of Grassmann-Cayley algebras), the right complement map (also known as Hodge duality), and the dual coordinate map (also known as \poincare duality). A more detailed comparison of Hodge and \poincare duality follows in Sec. \ref{sec:cph}. It is shown that in terms of a given basis, the underlying calculations are identical, differing only in how the results are interpreted: as elements of $G$ (Hodge) or of $G^*$ (\poincare). But $\PN$, in contrast to $\CR$, is a natural map. Finally, it introduces the terms single algebra, resp., double algebra representation to distinguish the results. Sec. \ref{sec:dngp} introduces \emph{geometric duality}, a powerful but little-known feature of the double algebra approach that plays a key role in euclidean rigid body mechanics.  Sec. \ref{sec:impdaap} then proposes a unified software implementation that offers advantages to both double and single algebra approaches.  Sec. \ref{sec:conc} sums up the results presented.

\section{Projective geometric algebra}
\label{sec:prel}
The discussion assumes basic familiarity with geometric algebra and linear algebra. 

Projective geometric algebra (\cite{gunn2011}, \cite{gunn2017a}, \cite{gunn2019}) introduces the distinction of \quot{point-based} or \quot{plane-based} geometric algebras, depending on whether 1-vectors are interpreted to be elements of $\RP{n}$ (points) or of its dual space $\RP{n}^*$ (hyperplanes). The wedge product $\wedge$ is accordingly \emph{join}, resp., \emph{meet}. These two approaches lead to the \emph{standard}, resp., \emph{dual} construction of space. In general we let $G$  represent any PGA -- either point- or plane-based -- and $G^*$ its dual algebra, representing the same projective space but using the opposite construction. In the context of PGA, then, every geometric primitive appears doubled: once in plane-based form and once in point-based form. The existence of this little-known but pervasive \emph{geometric duality} is explored more fully in Sec. \ref{sec:dngp}. It plays a key role in the considerations presented here.

This article focuses on the case of most practical interest: $ G = \pdclal{3}{0}{1}$, 3D euclidean PGA, which we also refer to as EPGA.  We assume that $G^*$ has the same signature and is hence the algebra  $\pclal{3}{0}{1}$, dual euclidean space (DES), discussed more in Sec. \ref{sec:des}. The metric of $G^*$ plays only a small role in this article since it is focused on non-metric computations. 



\section{Calculating coordinates}
\label{sec:cc}

As an introduction to the use of duality in geometric algebra, we begin by showing how the $P$ map above can be used to obtain the coordinates of 1-vectors with respect to a given basis. 
 Readers who are familiar with this topic are encouraged to skip ahead to the discussion of the right complement map in Sec. \ref{sec:rcm}.


\subsection{Calculating coordinates using reciprocal frames}

Let $G$ be an $n$-dimensional  geometric algebra such that $\eye^2\neq 0$.
Let  $\mb{e}_i$ be a basis of the 1-vectors $\bigwedge^1$. Then the reciprocal basis $\mb{e}^i$ of $\bigwedge^1$ is defined as: 
\begin{align}
    \mb{e}^i := (-1)^{(i-1)} (\mb{e}_1 \wedge \mb{e}_2 ... \wedge \widehat{\mb{e}_i} \wedge ... \wedge \mb{e}_n)\cdot \eye^{-1} \label{eq:rf}
\end{align}

Here $\widehat{\e{i}}$ means to omit the $\e{i}$ element from the wedge product, while $\cdot$ represents the generalized inner product of two blades, that is the lowest-grade part of their geometric product.

\textbf{Exercise:} $\e{i}$ is an orthonormal basis $\iff$ $\e{i} = \mb{e}^i$ for all $i$. 

One of the main purposes of the reciprocal frame is to quickly compute the coordinates of a vector with respect to a non-orthogonal basis (as discussed in Sec 3.8 of \cite{dfm07} and later in the discussion of directional derivative). 
In fact, the reciprocal basis is so constructed that if $\mb{b} = b_i \e{i}$ is a 1-vector, then $b_i = \mb{e}^i \cdot \mathbf{b} $. 
We now turn to a method to calculate coordinates that doesn't depend on the metric.

\subsection{Calculating coordinates without a metric}

There is no logical necessity for invoking a metric in order to calculate  coordinates with respect to a basis. They are well-defined in any real vector space. (\cite{gunn2017a}, Sec. 4.1.2). 

We present some results in the vector space setting; the transition to the projective setting doesn't present any difficulties since \quot{calculating coordinates} commutes with \quot{projectivization} (of course up to multiplication by a non-zero scalar!).

\subsubsection{Canonical dual basis}
\label{sec:cdb}
The key step is to construct the \emph{canonical dual basis}  in $\bigwedge^{n-1}$, the space of co-vectors. We use the notation $\{\db{i}\}$ to distinguish it from the reciprocal frame.  This is a basis $\{\db{i}\}$ for the $(n-1)$-vectors $\bigwedge^{n-1}$ such that $\e{i} \wedge \db{j} = \delta^j_i \eye$.   Define $\db{i} := (-1)^{(i-1)}\e{1} \wedge ... \wedge \mb{\widehat{e}_i} \wedge ... \e{n}$. 
Then it's easy to verify $\e{i} \wedge \db{j} = \delta^j_i\eye$ so it satisfies the definition of canonical dual basis.  

\myboldhead{Remark} This construction is identical to that of the reciprocal frame, but leaving out the pseudoscalar multiplication.

\myboldhead{The canonical dual basis is not coordinate-free}
To see this, in any 2D PGA, consider replacing  $\e{1}$ with $\e{0}+\e{1}$ so the new basis $\{\e{0}', \e{1}'\} = \{\e{0}, \e{0}+\e{1}\}$. Then $\mathbf{I}' = \mb{e}'_0 \wedge \mb{e}'_1 = \mathbf{I}$. Construct the canonical dual bases with respect to these two bases: for the first basis $\{\e{0}^*, \e{1}^*\} = \{\e{1}, -\e{0}\}$, while with respect to the second basis  $ \{\mathbf{e'}_0^*, \mathbf{e'}_1^*\} = \{\mathbf{e'}_1, -\mathbf{e'}_0\} = \{\e{0}+\e{1}, -\e{0}\}$ where the second equality follows from substituting. So  the dual of $\e{0} = \e{0}'$ with respect to the first basis is $\e{1}$ while with respect to the second basis it is $\e{1}' = \e{0}+\e{1}$. So this is not coordinate-free. 




\subsubsection{Calculating coordinates using the dual basis}

Let $\mb{b}$ be an arbitrary 1-vector $\mb{b} = \sum b_i\mb{e}_i$. That is, we know by general principles that the coordinates $b_i$ exist, the challenge is to calculate them. 

First, introduce the linear map $S: \bigwedge^n \rightarrow \mathbb{R}$ satisfying $\mb{p} = S(\mb{p})\eye$, the \emph{coordinate} of $\mb{p}$ with respect to the unit pseudoscalar $\eye$. Then define:

\[\hat{b}_i = S(\mb{b} \wedge \db{i}) ~~~\text{and} ~~~\widehat{\mb{b}} := \sum \hat{b}_i \mb{e}_i\] 

  We show that $\mb{b} = \widehat{\mb{b}}$.  Indeed, consider the $n$ quantities $(\mb{b}-\widehat{\mb{b}})\wedge \db{j}$.  $\mb{b}\wedge \db{j} = \hat{b}_j \eye$  by definition of $\hat{b}_j$, while $\widehat{\mb{b}} \wedge \db{j} = (\sum{\hat{b}_i \e{i}}) \wedge \db{j} = \hat{b}_j \eye$ by definition of dual basis. Hence $(\mb{b}-\widehat{\mb{b}})\wedge \db{j} = 0$ for all $j$, which can only happen if $\mb{b}-\widehat{\mb{b}}=0$, so $\mb{b} = \widehat{\mb{b}}$.

\myboldhead{Example} This metric-neutral approach faithfully mirrors the nature of coordinates as seen in this example. Begin with a non-orthogonal basis $\{e_i\}$ for $\R{3}$ and consider the dual basis, consisting of the 3 planes $\{e_k^* := e_i \wedge e_j\}$ spanned by pairs of these basis vectors.  Choose an arbitrary point $\mb{P}$.  Then $\mb{P}$ can be represented uniquely as the intersection of three planes obtained by parallel displacement of $\{e_k^*\}$. The coordinates of $\mb{P}$ are just the magnitude of these displacements. The wedge products $\mb{P}\wedge e_k^*$ that appear in the coordinate calculations measure these displacements.

\subsection{Comparison of the two methods}  Note that the metric is applied twice when using the reciprocal frame to calculate coordinates: first in calculating the reciprocal frame, and secondly when computing the individual coordinate components of the 1-vector. The approach using the dual basis basically omits the first inner product and replaces the second with a wedge product.  This is a typical pattern in the usage of pseudoscalar multiplication to compute non-metric quantities: by applying it twice, the metric dependence is nullified, since $\eye^2=\pm1$.

\subsection{Canonical dual basis for general $k$-vectors} 
\label{sec:cdbk}
The canonical dual basis can be defined for $\bigwedge^k$ for any $k$. Since this will be important for the sequel we introduce it here. We sketch the construction for general dimension $n$ but focus then on the case $n=4$.   The preferred condition remains that for a basis vector $e_I$ with index set $I \subset \{1,2,...,n\}$ (written as a $k$-digit integer) we want that $\e{I} \wedge \db{I} = \eye$ where $\db{I}$ is a basis vector of grade $(n-k)$. 

When $(n-k)>1$, we can calculate $\db{I}$ in the following way. Choose an index set $I^\perp$ such that the concatenation  $(II^\perp)$ is an even permutation of $(1234)$. For example when $I = 1$, $I^\perp = 032$. Since $\#(I^\perp) = (n-k) > 1$, we can always swap elements if necessary to obtain the desired condition. 
Then define $\db{I} := \e{I^\perp}$. For example, when $I = \{02\}$, then $I^\perp = 31$ since $(0231)$ is an even permutation, and $\db{I^\perp} = \e{31}$. 

When $k = n-1$, $\#(I^\perp) = 1$, say $I^\perp = \{i\}$. Then  define $\db{I} := (-1)^{(n-1)}\e{i}$. 

Consult the row of Table \ref{tab:j3d} labeled $\CR$ for a full listing for the case $n=4$. 

\myboldhead{Correction} The original description of the calculation of this index set in \cite{gunnthesis}, Sec. 2.3.1.3, overlooked the special case $k = n-1$ and the consequent introduction of some minus signs in the dual basis vectors.

\section{The right complement  map $\CR$}
\label{sec:rcm}

In the context of a geometric algebra the canonical dual basis is closely related to the \emph{right complement} map, which plays a key role in the discussion of the regressive product below. The same results can also be obtained using \textit{left} complements; the two approaches are equivalent. We turn now to discuss the right complement map.


A \quot{right complement} of an element $\mathbf{x}$ in a geometric algebra $G$ is an element 
$\mathbf{x}^*$ satisfying $\mathbf{x} \mathbf{x}^* = \mathbf{I}$. When $\mathbf{x}^{-1}$ exists, then this right complement is unique and can be directly computed as $\mathbf{x}^* = \mathbf{x}^{-1}\mathbf{I}$. 
However, when the metric is degenerate, the right complement is no longer well-defined.

\myboldhead{Projective formulation} In projective geometry, two elements $\mb{x}$ and $\mb{y}$ are considered equivalent when $\mb{x} = \lambda \mb{y}$ for $\lambda \in \mathbb{R}, \lambda \neq 0$. Since the focus of this article is projective geometric algebra, we rephrase the condition of right complement to be $\mathbf{x} \mathbf{x}^* = \lambda \mathbf{I}$ for some $\lambda \neq 0$. Of course it's always possible
to scale $\mb{x}^*$ by $\lambda^{-1}$ so that $\mb{x}\mb{x}^* = \eye$. In this case we say $\mb{x}^*$ is the \emph{strict} right complement of $\mb{x}$. Unless we say otherwise we assume in the sequel we are dealing with the projective right complement.

\myboldhead{Example} 
When the metric is degenerate, as with euclidean PGA, and $\mathbf{x}^2=0$, then the right complement is not well-defined. For example, let $\mathbf{x} = -\e{021}$, the ideal point in the negative z-direction. Then $\mathbf{x}^* := \e{3}$ ($\e{3}$ is the z-plane) certainly satisfies $\mathbf{x}\mathbf{x}^*=\mathbf{I}$ but so does $\mathbf{x}_\alpha^*:=\e{3} + \alpha \e{0}$ for any $\alpha \in \mathbb{R}$ since $\e{0}^2=0$. 

\begin{figure}
  \centering
  \def\ddc{0.6}
  \def\ddd{0.35}
    \includegraphics[ width=\ddc\columnwidth] {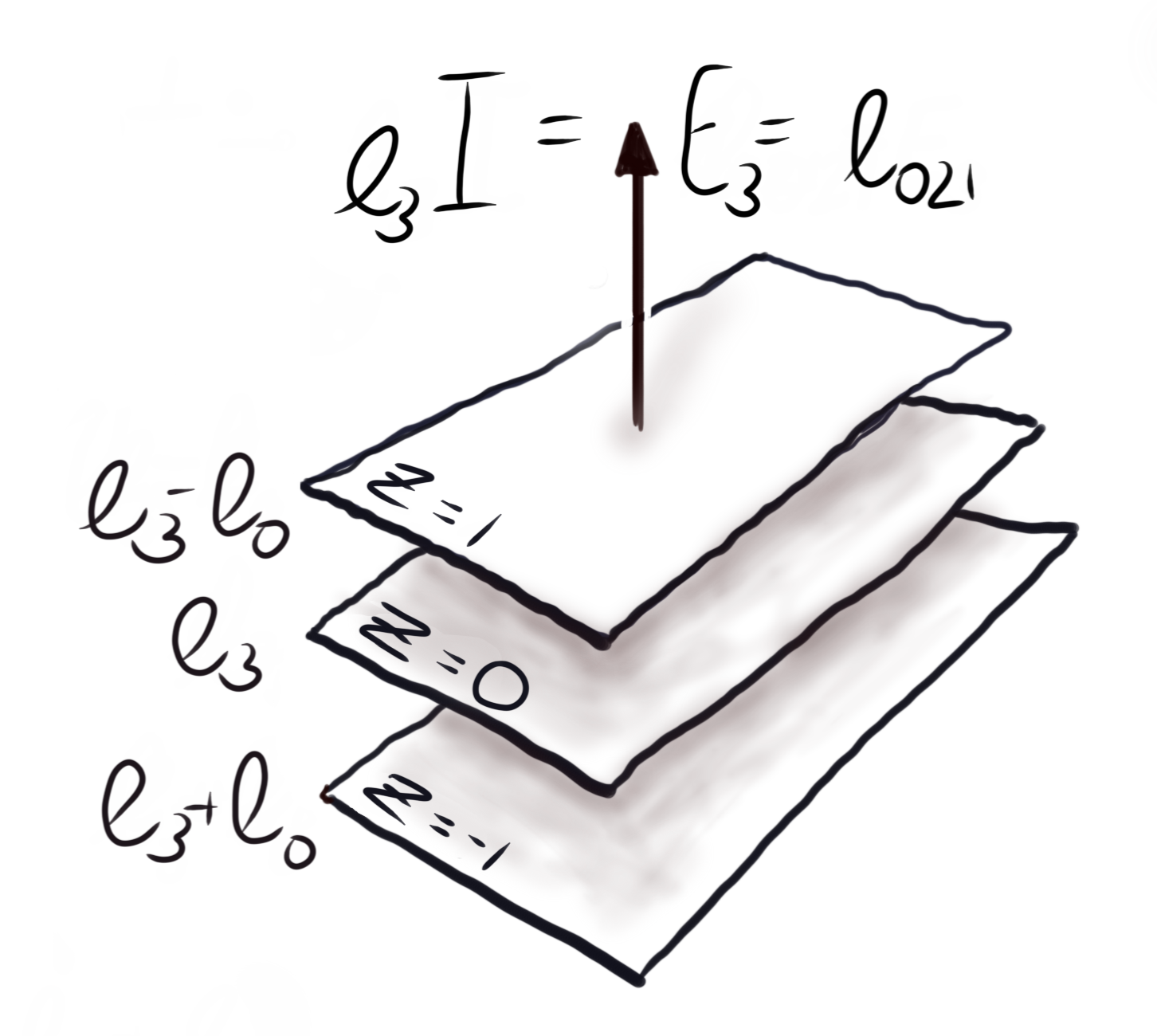} 
  \caption{Right complement  of $-\e{021}$ is not well-defined in $\pdclal{3}{0}{1}$. All the planes of the form $z=\alpha$ are right complements. }
\label{fig:parplanes}
\end{figure}
\myboldhead{Geometric interpretation} $\mb{x}_\alpha^*$ represents the plane parallel to $\e{3}$ at a signed distance of $\alpha$ from it.  We conclude that if a plane $\mb{p}$ is a right complement then all planes parallel to $\mb{p}$ are also a right complement.  Hence \textbf{the} \quot{right complement} isn't well-defined in this case. See Fig. \ref{fig:parplanes}.

\subsection{Defining \textbf{the} right complement using the dual basis}
\label{sec:drcu}
If we want to have a well-defined right complement map in the presence of null elements it's clear from the above that it can't depend on the metric. In the following we obtain a non-metric form of right complement. The price, as we shall see, is that it is not coordinate-free.

Define the \emph{right complement} map $\CR$ \textit{with respect to} the basis $\e{i}$ as follows. If $
\mathbf{b} = \sum b_i \e{i}$, then 
\[\CR: \bigwedge^{1} \rightarrow \bigwedge^{n-1}: \CR(\mathbf{b}) := \sum{b_i \db{i} }\]
where $\db{i}$ is the canonical dual basis defined in Sec. \ref{sec:cdb}.

Verify that 1) $\mb{b} \wedge H(\mb{b}) = \lambda \eye$, where $\lambda := \sum{b_i^2}$, hence 2) $\lambda^{-1}H(\mb{b})$ is the strict right complement, and 3) $H(\mb{b})$ is the projective right complement of $\mb{b}$.

\subsubsection{Geometric interpretation}
Notice that this converts a 1-vector to an $(n-1)$-vector by interpreting its coordinates with respect to the canonical dual basis of $(n-1)$-vectors. Or, expressed in other words, it is the linear extension of the dual basis map $\e{i} \rightarrow \db{i}$ to the whole vector space.
Since it is defined in terms of the canonical dual basis, $\CR$ is  not coordinate-free. 

\subsubsection{Right complement and Hodge duality} 
\label{sec:rcahd}
It can be shown that $\CR$ induces a non-degenerate inner product $\langle , \rangle_H$ on the 1-vectors such that \[ \mbf{X}  \wedge \CR(\mbf{Y}) = \langle \mbf{X}, \mbf{Y}\rangle_H \eye\] In this form it is well-known as the \emph{Hodge star} map, or \emph{Hodge duality}.  This in effect overcomes the absence of a non-degenerate metric by creating a new one. We'll see below in Sec. \ref{sec:pd} how \poincare duality avoids this and is also coordinate-free.
The extension of $\CR$ to the full algebra is not difficult and can be found in standard texts. 

\section{Calculating the regressive product}
\label{sec:regpro}
We assume that $G$ is euclidean PGA, that is, the plane-based algebra $G = \pdclal{3}{0}{1}$, so that $\wedge$ is meet. We want to investigate how to compute the join of two elements. To distinguish this operation from the standard wedge operation $\wedge$ which is sometimes called the \emph{progressive} product, it is called the \emph{regressive} product and written $\vee$.  (Similar remarks apply in a standard PGA; there the 1-vectors are points and the regressive product is join, but since euclidean PGA is the focus of this article, we focus on a plane-based $G$).

Notice that the regressive product $\mb{X} \vee \mb{Y}$ of a $k$- and an $m$-vector will be 0 unless $k+m \geq 4$.

We begin by describing how the regressive product has been traditionally computed using pseudoscalar multiplication, followed by descriptions of three non-metric approaches, all of which yield the same results.

\subsection{Via non-degenerate metric}
\label{sec:ndm}

  In the terminology of projective geometry, the map $P$ is the \emph{\textbf{p}olarity on the metric quadric}. As such, it maps a subspace onto its orthogonal complement.  
When $\eye^2\neq0$, then $P^2(X) \cong X$ (up to projective equivalence) and maps each point to a unique plane and vice-versa, that is, it is a grade-reversing involution. It has been traditionally used to calculate the regressive product in either of two forms:

\[\mb{X}\vee\mb{Y} := \mb{I}^{-1}(\mb{I}\mb{Y} \wedge \mb{I}\mb{X})
~~~~~(= (\mb{Y} \mb{I}^{-1} \wedge \mb{X} \mb{I}^{-1})\mb{I})\]

This formula is widely available in the literature, e.g., Formula 21.8 in \cite{dorst2011}, p. 443.

  However, when the metric is degenerate, $\mb{I}^{-1}$ doesn't exist, and $P$ is not injective. In EPGA, for example, $P$ collapses the whole of euclidean space onto the ideal plane: the orthogonal complement of any euclidean multi-vector is ideal.  A plane is mapped to its normal direction (an ideal point); a line is mapped to the ideal line of directions perpendicular to the direction of the line; every point is mapped to the ideal plane. This behavior is related to that discussed in Sec. \ref{sec:rcm} and shown in Fig. \ref{fig:parplanes}. In this and other degenerate cases non-metric methods are required to obtain the regressive product.

\subsection{Via shuffle product}

Historically, Grassmann-Cayley algebras were developed  to meet this need in a coordinate-free way through the definition of a so-called \emph{shuffle} product (\cite{selig05}, Ch. 10, \cite{white1995}). For an $j$- and $k$-blades $\mb{a} = a_1 \wedge a_2 \wedge ... \wedge a_j$ and $\mb{b} = b_1 \wedge b_2 \wedge ... \wedge b_k$ where $j+k \geq n$, the shuffle product is defined as:
\begin{multline*} 
\mb{a} \vee \mb{b} := \dfrac{1}{(n-k)!(j+k-n)!} \sum_\sigma \text{sign}(\sigma) (\text{det}(a_{\sigma(1)}, ..., a_{\sigma(n-k)}, b_1, ... b_k)a_{\sigma(n-k+1)}, ... , a_{\sigma(j)} 
\end{multline*}
where $\sigma$ is a permutation of ${1,2,...,j}$.

\myboldhead{Exercises} Calculate the shuffle products  a) $\e{123} \vee \e{032}$ and b) $\e{0} \vee \e{123}$. [Answers: a) $\e{23}$, b) 1].

\subsection{Via dual coordinate map (\poincare duality)}
\label{sec:pd}

\cite{gunnthesis}, \cite{gunn2019} describes an alternative duality map known in the world of multilinear algebra as \poincare duality $\PN: G \rightarrow G^*$ (\cite{greub2-67}).  The action of  $\PN$ is the essence of simplicity: it maps an element in the first algebra to the element of the dual algebra that represents the \textbf{same} geometric entity.  Since it converts back and forth between coordinates in the point-based and plane-based algebras, $\PN$ has historically been known also as the \emph{dual coordinate map}.  Like $P$, $\PN$ is a grade-reversing map that is a vector space isomorphism on each grade. 
We'll see below in Sec. \ref{sec:natural} that $\PN$ is essentially an \textit{identity map} and plays an essential role in the discussion of geometric duality in Sec. \ref{sec:dngp}.

To better understand $\PN$ we look more closely at its action in projective space $\RP{3}$. 

\subsubsection{The fundamental tetrahedron}
\begin{figure}
  \centering
  \def\ddc{0.8}
  \def\ddd{0.48}
    \setlength\fboxsep{0pt}\fbox{\includegraphics[width=\ddc\columnwidth]{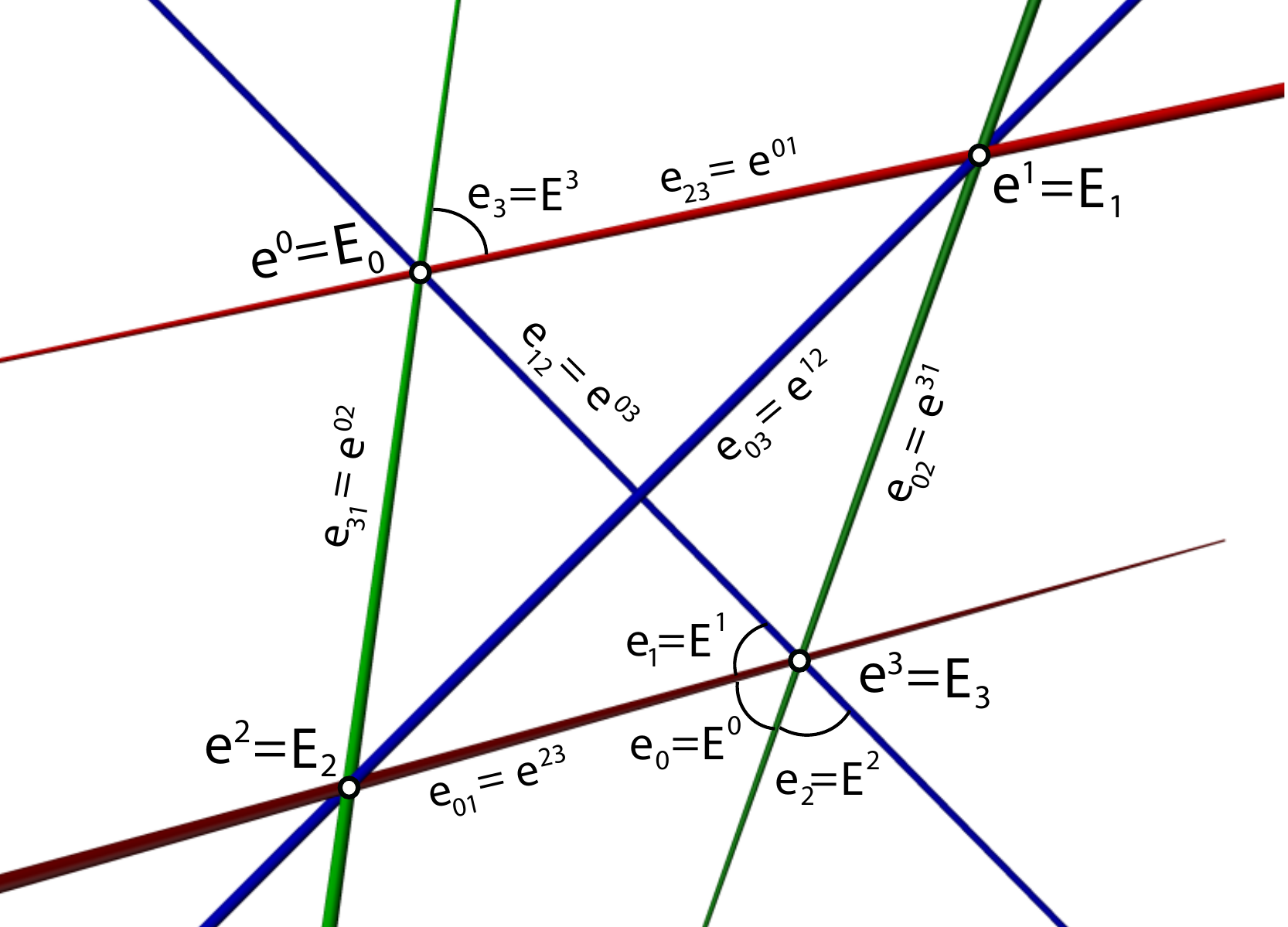}} 
  \caption{The fundamental tetrahedron for $\RP{3}$.}
\label{fig:fundtet}
\end{figure}

$\PN$ is defined by its action on the standard bases shown in 
Fig. \ref{fig:fundtet}.
This shows the \textit{fundamental tetrahedron} of projective 3-space $\RP{3}$.  The vertices, edge-lines, and face-planes of the fundamental tetrahedron represent the basis vectors for the exterior algebra.  Superscripts are used for the standard construction, where  the vertices are 1-vectors, edge-lines are 2-vectors, and face-planes are 3-vectors, while subscripts are used for the dual construction where planes are 1-vectors, lines are 2-vectors, and points are 3-vectors. 

3-vectors are represented by capital letters according to the rule $\mbf{E}_i := \mbf{e}_{jkl}$, $\mbf{E}^i:=\mbf{e}^{jkl}$ whereby the defining condition is that $(ijkl)$ is an even permutation of $(0123)$ and is alphabetically the least such permutation: $\mbf{E}_1 = \mbf{e}_{032}$. Two-vectors are paired in this way with 2-vectors. For example, the left green line is $\e{31} = \e{3} \wedge \e{1}$ in $G$ (plane-based coordinates, $\wedge$ is meet) while in the dual algebra the same line is represented by $\mbf{e}^{02}$ since $\mbf{e}^{02} = \mbf{e}^0 \wedge \mbf{e}^2$ in $G^*$ (point-based coordinates, $\wedge$ is join). We use $\e{31}$ since $(0231)$ is an even permutation. Table \ref{tab:j3d} shows a complete listing of $\PN: G \rightarrow G^*$. Note that the index sets that appear are exactly the same as those produced by $P$, the (non-deganerate) metric polarity. 

In EPGA, the plane $\e{0}$ is the ideal plane; the point $\mb{E}_{0}$ is the origin, the lines $\mb{e}^{0i}$ for $i\in\{1,2,3\}$ are the coordinate axes, etc.

\begin{table}[b]
\centering

\colorbox{ltpink}{
\resizebox{\textwidth}{!}{
\begin{tabular}{|l | c|c|c|c|c|c|c|c|c|c|c | c | c | c | c | c |}
\toprule
 $\mbf{e}$ & $\mathbf{1}$  & $ \e{0}$ & $\e{1}$ & $\e{2} $ & $\e{3} $ & $\e{01} $ & $\e{02} $ & $\e{03}$ & $\e{12} $ & $\e{31} $ & $\e{23} $ & $\e{123} $ & $ \e{032} $ & $ \e{013} $ & $ \e{021} $ & $ \eye $  \\
\midrule
\rowcolor{ltred}
$\MP(\mbf{e})$ & $ \e{0123} $ & $ \e{123}$ & $ \e{032} $ & $\e{013} $ & $\e{021} $ & $\e{23} $ & $\e{31} $ & $\e{12}$ & $\e{03} $ & $\e{02} $ & $\e{01} $ & $-\e{0} $ & $ -\e{1} $ & $ -\e{2} $ & $ -\e{3} $ & $ \mathbf{1} $ \\
\midrule
\rowcolor{ltyellow}
$\PN(\mbf{e})$ & $ \ej{0123} $ & $ \ej{123}$ & $ \ej{032} $ & $\ej{013} $ & $\ej{021} $ & $\ej{23} $ & $\ej{31} $ & $\ej{12}$ & $\ej{03} $ & $\ej{02} $ & $\ej{01} $ & $-\ej{0} $ & $ -\ej{1} $ & $ -\ej{2} $ & $ -\ej{3} $ & $ \mathbf{1} $\\
\midrule
\rowcolor{ltgreen}
$\CR(\mbf{e})$ & $ \e{0123} $ & $ \e{123}$ & $ \e{032} $ & $\e{013} $ & $\e{021} $ & $\e{23} $ & $\e{31} $ & $\e{12}$ & $\e{03} $ & $\e{02} $ & $\e{01} $ & $-\e{0} $ & $ -\e{1} $ & $ -\e{2} $ & $ -\e{3} $ & $ \mathbf{1} $\\
\bottomrule
\end{tabular}
}
}
\caption{Metric polarity $\MP$, dual coordinate $\PN$ and right complement $\CR$ maps for $\RP{3}$.}
\label{tab:j3d}
\end{table}

\subsubsection{$\PN$ in terms of a basis}
\label{sec:pnb}
Expressing $\PN$ in terms of a basis proceeds exactly as the definition of $\CR$ described in Sec. \ref{sec:cdbk}. However where $\CR(\e{I}) := \e{I^\perp}$, one has $\PN(\e{I}) := \mb{e}^{I^\perp}$: the target domain is $G^*$, not $G$.   See Tab. \ref{tab:j3d}, the row labeled $\PN$. 

Consider $\e{12}$, the meeting plane of $\e{1}$ (the $x=0$ plane) and $\e{2}$ (the $y=0$ plane).  According to the table, $\PN(\e{12}) = \mb{e}^{03}$, the joining line of $\mb{e}^0$, the origin, and $\mb{e}^3$, the ideal point in the $z$-direction. Our choice of $I^\perp$ as complementary to $I$ guarantees that both points are contained in both planes $\e{1}$ and $\e{2}$ and hence join to form their common line $\e{12}$. 

In general, if $\e{I}$ represents the meet of the $k$ planes $\mathfrak{P} := \{\e{i_1}, \e{i_2}, ... \e{i_k}\}$, then the $(n-k)$ points represented by the indices of $I^\perp$ will be exactly the basis vectors of $G^*$ that are incident with all of the planes in $\mathfrak{P}$. Hence their join lies in $\e{I}$. Reversing the argument establishes that the meet must lie in the join, so that $\e{I}$ and $\PN(\e{I})$ are the same geometric entity.


\subsection{Via right complement map (Hodge duality)}
The right complement map $\CR$ or Hodge duality, provides another  way to calculate the regressive product when the metric of $G$ is degenerate.   Again, the formula is isomorphic to the preceding two:
\[\mb{X}\vee\mb{Y} := {\CR}^{-1}(\CR(\mb{X}) \wedge \CR(\mb{Y}))\]

As noted above in Sec. \ref{sec:rcahd}, Hodge duality is equivalent to defining a second (non-degenerate) metric on $G$.
 We discuss that and other related issues in the next section.

\subsection{Comparison of the different approaches}
\label{sec:compare}
The shuffle product is quite different from the other approaches since it calculates the regressive product directly from its arguments. Leaving it to the side,  the remaining three approaches can be expressed using some progressive product:
\begin{align*}
\text{Metric polarity:}~~~ &\mb{X}\vee\mb{Y} := \mb{I}^{-1}(\mb{I}\mb{Y} \wedge \mb{I}\mb{X}) \\
\text{\poincare duality:}~~~ &\mb{X}\vee\mb{Y} := {\PN}^{-1}(\PN(\mb{X}) \wedge \PN(\mb{Y})) \\
\text{Hodge duality:}~~~ &\mb{X}\vee\mb{Y} := {\CR}^{-1}(\CR(\mb{X}) \wedge \CR(\mb{Y}))
\end{align*}
Notice however that the first formula reverses the order of its arguments. This  reflects the fact that the metric polarity is expressed using the geometric product while the other two forms of duality involve algebra maps. $\PN$ and $\CR$ differ in their target algebra.
Table \ref{tab:j3d} contains a full listing of the action of all three duality maps on a basis set. The index sets for the image basis set agree in all cases. 

\begin{table}[b]
\centering

\colorbox{ltpink}{
\resizebox{\textwidth}{!}{
\begin{tabular}{|l | c|c|c|c|c|c|c|c|c|c|c | c | c | c | c | c |}
\toprule
 $\mbf{e}$ & $\mathbf{1}$  & $ \e{0}$ & $\e{1}$ & $\e{2} $ & $\e{3} $ & $\e{01} $ & $\e{02} $ & $\e{03}$ & $\e{12} $ & $\e{31} $ & $\e{23} $ & $\e{123} $ & $ \e{032} $ & $ \e{013} $ & $ \e{021} $ & $ \eye $  \\
\midrule
\rowcolor{ltred}
$\MP^{-1}(\mbf{e})$ & $ \e{0123} $ & $ -\e{123}$ & $ -\e{032} $ & $-\e{013} $ & $-\e{021} $ & $\e{23} $ & $\e{31} $ & $\e{12}$ & $\e{03} $ & $\e{02} $ & $\e{01} $ & $\e{0} $ & $ \e{1} $ & $ \e{2} $ & $ \e{3} $ & $ \mathbf{1} $ \\
\midrule
\rowcolor{ltyellow}
$\PN^{-1}(\mbf{e})$ & $ \ej{0123} $ & $ -\ej{123}$ & $ -\ej{032} $ & $-\ej{013} $ & $-\ej{021} $ & $\ej{23} $ & $\ej{31} $ & $\ej{12}$ & $\ej{03} $ & $\ej{02} $ & $\ej{01} $ & $\ej{0} $ & $ \ej{1} $ & $ \ej{2} $ & $ \ej{3} $ & $ \mathbf{1} $\\
\midrule
\rowcolor{ltgreen}
$\CR^{-1}(\mbf{e})$ & $ \e{0123} $ & $ -\e{123}$ & $ -\e{032} $ & $-\e{013} $ & $-\e{021} $ & $\e{23} $ & $\e{31} $ & $\e{12}$ & $\e{03} $ & $\e{02} $ & $\e{01} $ & $\e{0} $ & $ \e{1} $ & $ \e{2} $ & $ \e{3} $ & $ \mathbf{1} $\\
\bottomrule
\end{tabular}
}
}
\caption{Inverses of the $\MP$, $\PN$ and $\CR$ maps for $\RP{3}$.}
\label{tab:j3di}
\end{table}

For future reference the inverse maps $\MP^{-1}$, ${\PN}^{-1}$, and ${\CR}^{-1}$ are given in Table \ref{tab:j3di}. Note the slight differences to Table \ref{tab:j3d} occasioned by the minus signs in the latter.

\section{Comparison of \poincare and Hodge duality}
\label{sec:cph}
Since we are interested here with non-metric duality, we now compare the duality maps $\PN$ and $\CR$. 
Tab. \ref{tab:j3d} reveals that, expressed in coordinates, the two maps are identical but have different target spaces. Accordingly, the only difference is that the index set of the values appears in $\PN$ as \textit{superscripts} while $\CR$ uses them as \textit{subscripts}. 

In light of this, it might seem surprising that only $\PN$ is coordinate-free. How can we understand this? We define a map $\RL: G \rightarrow G^*$ by $\RL = \PN \circ \CR^{-1}$. Consult the diagram below 
In the special coordinates that we have introduced above, $\RL$ simply \textit{\textbf{r}aises the indices}. What is $\RL$ doing geometrically? 

\begin{center}
\begin{tikzcd}
\bigwedge^{n-k}\subset G\arrow[r, "R"] & (\bigwedge^{n-k})^*\subset G^* \\
\bigwedge^k \subset G \arrow[u, "H"] \arrow[ur,"J"] & 
\end{tikzcd}

\end{center}

\begin{theorem}
$\RL(\mbf{X})$ is the left complement of $\PN(\mbf{X})$ in $G^*$.
\end{theorem}
\begin{proof}
We write the wedge product in $G^*$ as $\wedge_*$ and the pseudoscalar of $G^*$ as $\eye_*$.
\begin{align}
\RL(\mbf{X}) \wedge_* \PN(\mbf{X}) &= \\
\RL(\mbf{X}) \wedge_* \RL(\CR(\mbf{X})) &= \\
\RL(\mbf{X} \wedge \CR(\mbf{X})) &= \RL(\eye) = \eye_*
\end{align}
The step from (3) $\rightarrow$ (4) follows from the fact that $\RL$ is an algebra isomorphism. \end{proof}\textbf{Example:} When $\mbf{X} = \e{1}$, then $\RL(\mbf{X}) = \PN(\CR^{-1}(\mbf{X})) = \PN(-\e{032}) = \ej{1}$,  $\PN(X) = \e{032}$, and $\RL(\mbf{X}) \wedge_* \PN(\mbf{X}) = \eye_*$.

We present without proof the following restatement of the theorem:
\begin{corollary}
$J^{-1}(R(\mb{X}))$ is the left complement of $\mb{X}$ in G.
\end{corollary}
That is, the $J$ pull-back of $R$ is the inverse of $H$, hence the left-complement in $G$.


 \subsection{$\PN$ is a natural map}
 \label{sec:natural}
Writing $\PN = \RL \circ \CR$, we can see that the above theorem provides another proof that $\PN$ is the \quot{identity map} on geometric primitives, since, at the level of geometry, it is the composition of right complement with left complement, which is trivially the identity map. More importantly this factorization makes clear why $\PN$ is coordinate-free and $\CR$ is not: whatever coordinate dependence is present in $\CR$ (right complement) is undone by the application of $\RL$ (left complement).  Or, in other words, $\PN$ is a \emph{natural} map while $\CR$ is not.



\subsection{Double algebra versus single algebra}
We can characterize \poincare duality as a \emph{double algebra} approach since every invocation of the dual map $J(X)$ moves the geometric primitive $X$ back and forth between $G$ and $G^*$.  Invocations of $H(X)$, Hodge duality, on the other hand, remains in the algebra $G$ and is equivalent, as noted above, to multiplication by a second, non-degenerate pseudoscalar. 
It turns out that the double algebra approach is equivalent to   \emph{geometric duality}, a little-known but powerful feature of projective geometry that we now turn to discuss.

\section{Geometry duality: the double nature of geometric primitives in PGA}
\label{sec:dngp}
\begin{figure}
  \centering
  \def\scl{0.8}
  \def\ddc{0.348}
  \def\ddd{0.384}
    \setlength\fboxsep{0pt}\fbox{\includegraphics[width=\ddd\columnwidth]{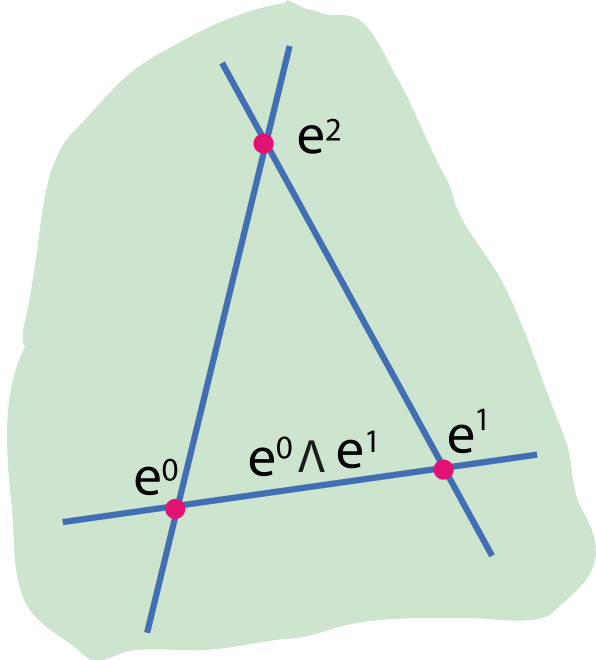}} \hspace{0.02\columnwidth}
    \setlength\fboxsep{0pt}\fbox{\includegraphics[width=\ddc\columnwidth]{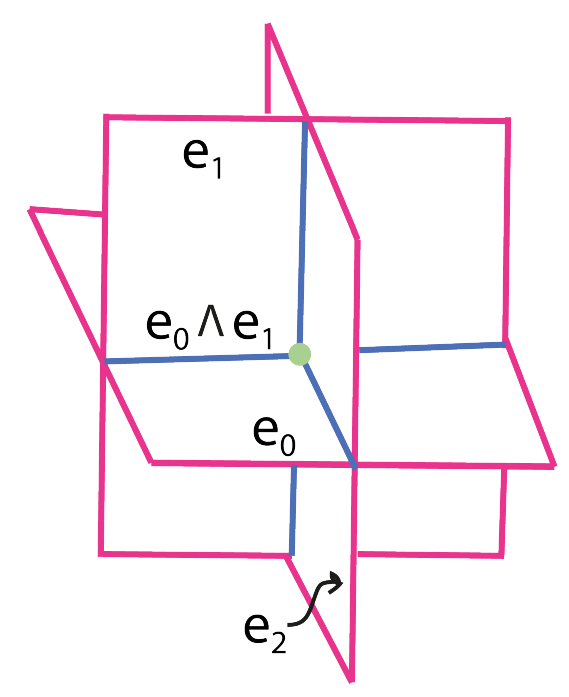}} 
  \caption{Building up space with the wedge operator in point-based (\emph{left}) and plane-based (\emph{right}) algebras.}
\label{fig:dualaufbau}
\end{figure}

As everyone who has begun to learn it has found out, one of the challenges of learning euclidean PGA $G = \pdclal{3}{0}{1}$ is that it is a plane-based algebra. See Fig. \ref{fig:dualaufbau}.
1-vectors represent planes, the wedge product ($\wedge$) implements the \textbf{meet} operator and 
the regressive product $\vee$ is  the \textbf{join} operator. That means that a line $\ell \in \bigwedge^2 \subset G$ arises as the intersection of two planes, and that a point $\mb{P} \in \bigwedge^3 \subset G$ arises as the intersection of three planes. The plane, as building block, is an indivisible, 0-dimensional primitive.

We are much more used to considering 1-vectors as points, and constructing lines and planes by joining points together.  Indeed, this approach leads to the common assumption that points are 0-dimensional, lines are 1-dimensional, and planes are 2-dimensional. From this perspective, it might appear that in a plane-based algebra, we are building up objects of lower dimension from those of higher dimension, for example, a line is the meet of two planes.  How can we understand this paradox?

A clear definition of \emph{dimension} resolves this dilemma. We turn now to that task. 

\subsection{Point-based lines and planes} For example the join $m := P_1 \wedge P_2$ of two points in a point-based algebra is a 2-vector representing the joining line of the two points. The line is 1-dimensional since an arbitrary point $P $ incident with $m$ (that is, satisfying $m \wedge P = 0$) can be represented by the family of points $P_1 + \beta P_2$ for $\beta \in \mathbb{R}$ (to include $P_2$ in this expression you can use 1-dimensional homogeneous coordinate $\alpha : \beta$ in $\alpha P_1 + \beta P_2$). The line $m$ in this sense is conceived of as consisting of all the points incident with it, called a \emph{point range} in classical projective geometry. It's also an example of an  \emph{outer product null space} (OPNS) in the terminology of exterior algebra. Similar remarks apply to a plane $p = P_1 \wedge P_2 \wedge P_3$ as a 2-dimensional set of points, called a \emph{point field} in projective geometry. To sum up: the dimension $d$ of a geometric primitive in the standard construction depends on  how many linearly independent points $n$ you need to generate the primitive ( $d=n-1$).

\begin{table}
\centering

\colorbox{bkgd2}
{
\begin{tabular}{|| c ||c | c || c | c ||}
\toprule
 \textbf{PRIMITIVES} & \multicolumn{2}{c||} {Point-based}  & \multicolumn{2}{c||} {Plane-based} \\
\toprule
generic name & name & dimension & name & dimension \\
\midrule
point & \cellcolor{v1color} point & \cellcolor{v1color} 0  & \cellcolor{v3color} plane bundle & \cellcolor{v3color} 2 \\
\midrule
line &  \cellcolor{v2color} point range (spear) & \cellcolor{v2color} 1  & \cellcolor{v2color} plane pencil (axis) & \cellcolor{v2color} 1 \\
\midrule

plane & \cellcolor{v3color} point field & \cellcolor{v3color} 2  & \cellcolor{v1color} plane & \cellcolor{v1color} 0 \\
\midrule
\bottomrule
\end{tabular}
}
\caption{Double nature of geometric primitives in 3D PGA.}
\label{tab:doubleg}
\end{table}

\subsection{Plane-based lines and points} What does this look like in a plane-based algebra? Now we build up all the other primitives out of planes.  A plane in the dual construction is just as simple, indivisible, and 0-dimensional as a point is in the standard construction. Planes are combined by intersecting them.  Let $m = p_1 \wedge p_2$ be the intersection line of two planes. A third plane $p$ is incident with this line if it satisfies  $m \wedge p = 0$.  This set of incident planes is called a \emph{plane pencil} and is exactly analogous to the point range defined above. Just as in the point-based algebra we say that \quot{the joining line consists of all the points incident with it} in the plane-based algebra we say that \quot{the intersection line consists of all the planes incident with it.}  This pencil can be parametrized (exactly as above) as $\alpha p_1 + \beta p_2$ for $\{\alpha, \beta\} \in \RP{1}$, so it's 1-dimensional (with respect to the 0-dimensional building block, the plane).  

Similar remarks apply to the wedge of 3 planes to produce a point: $P = p_1 \wedge p_2 \wedge p_3$. The set of planes $p$ incident to this point satisfy $p \wedge P =  0$.  Just as the plane in the standard construction can be thought of as all the points incident to it, the point in the dual construction can be thought of as consisting of all the planes that are incident with it.  This set of planes is dual to the \emph{point field} above and is called a \emph{plane bundle} in projective geometry.

\begin{figure}[h]
  \centering
  \def\ddc{0.35}
  \def\ddd{0.48}
    \setlength\fboxsep{0pt}\fbox{\includegraphics[width=\ddc\columnwidth]{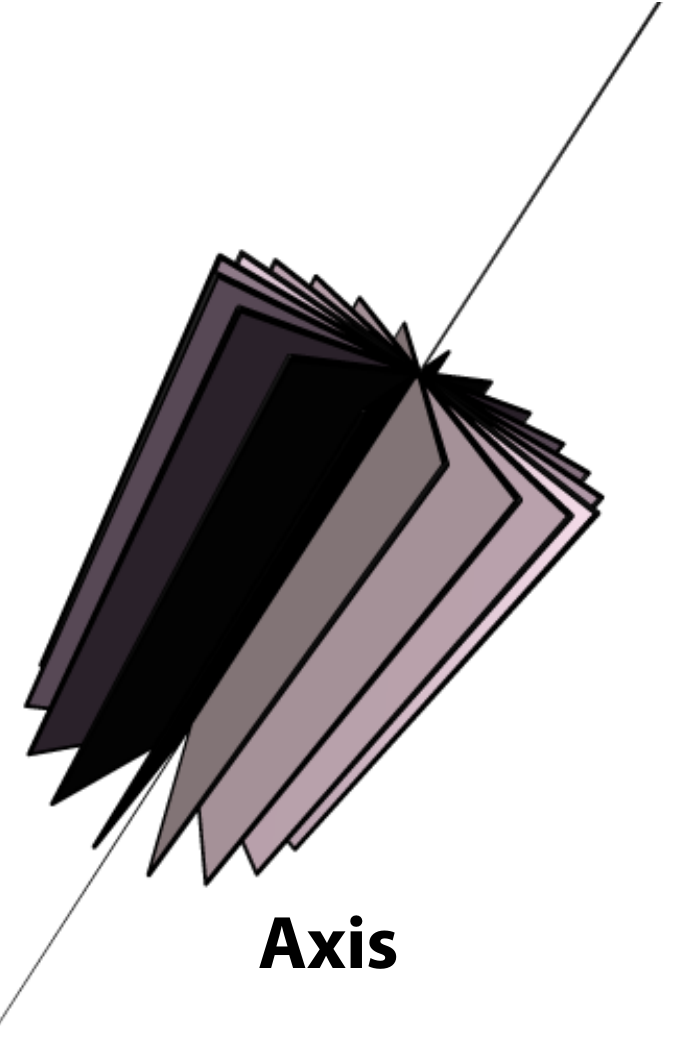} \hspace{0.1\columnwidth} \includegraphics[width=\ddc\columnwidth]{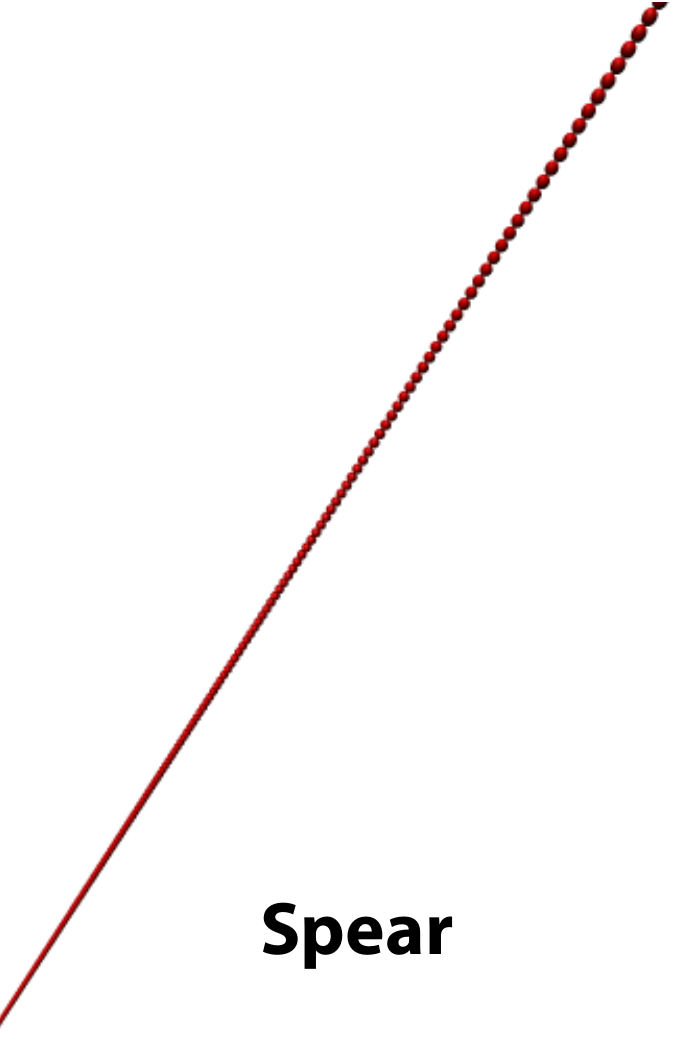}}
  \caption{Line as axis (plane pencil) and as spear (point range).}
\label{fig:spearaxis}
\end{figure}

Fig. \ref{fig:spearaxis} shows the difference between a line (bi-vector) in the plane-based algebra (where it is a plane pencil, or \emph{axis}) and its representation in a point-based algebra (as a point range, or \emph{spear}). We'll come back to this below when we see that this distinction plays a key role in rigid body mechanics.



\begin{figure}[ht]
  \centering
  \def\ddc{0.9}
  \def\ddd{0.48}
 \setlength\fboxsep{0pt}\fbox{\includegraphics[width=\ddc\columnwidth]{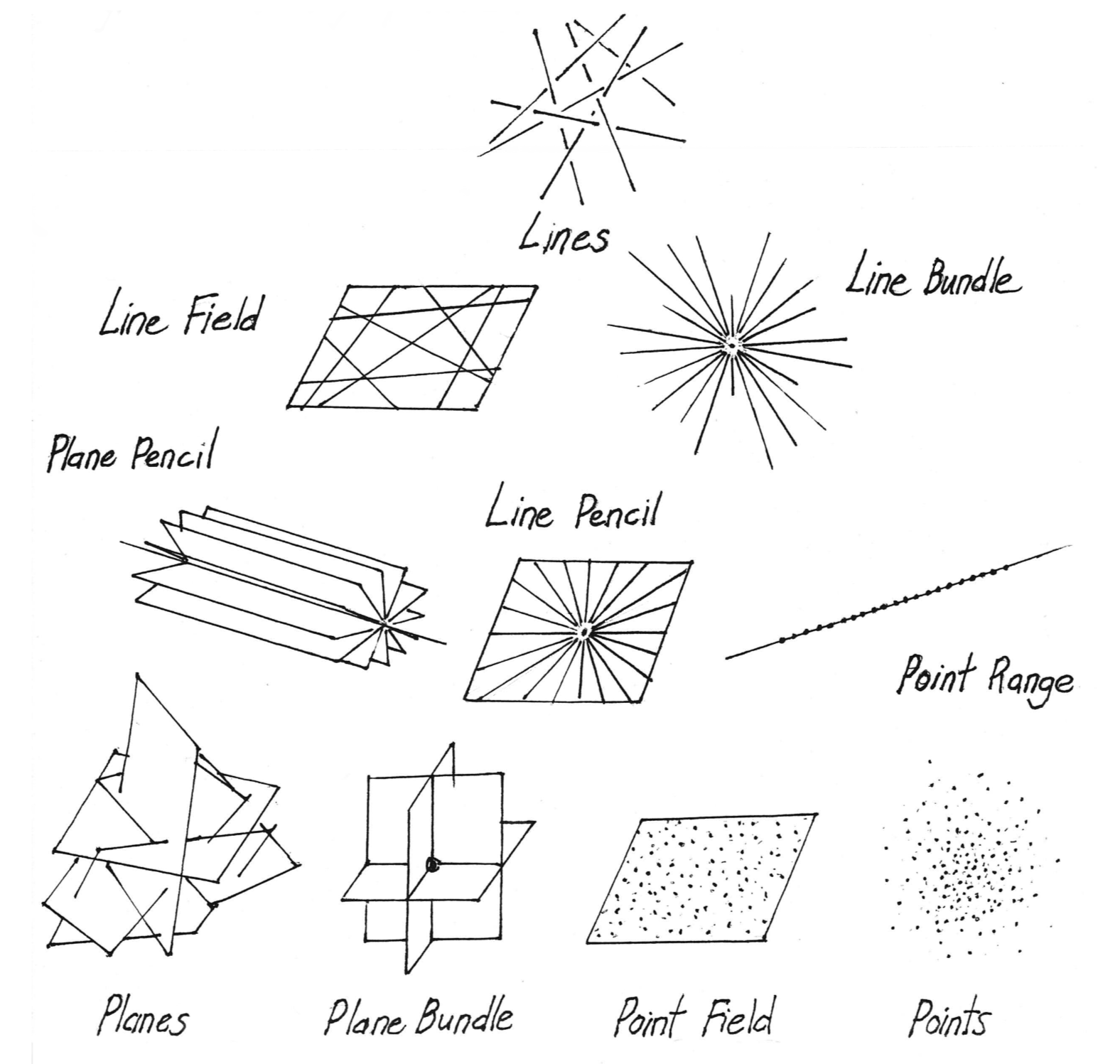}}
  \caption{The elemental triangle showing the 3 spaces (corners), the 6 composite primitives (sides), and the line pencil (middle).}
\label{fig:eltri}
\end{figure}

\subsection{There are no lines in 3D PGA}
The above line of thought leads to the insight that in geometric algebra, all geometric primitives appear in a double form, once built out of points, once out of planes. See Table \ref{tab:doubleg} and Figures \ref{fig:dualaufbau} and \ref{fig:spearaxis}. In a point-based algebra, the point \emph{qua} point is present, while in a plane-based algebra, the plane \emph{qua} plane is present.  But in neither algebra is a line \emph{qua} line present; every line is a one-dimensional composite object; in a point-based algebra it's a spear (point range), in a plane-based algebra it's an axis (plane pencil).  The dual coordinate map $\PN$ of PGA is exactly what is needed to switch easily back and forth between these two representations.


\subsection{The elemental triangle and the 10 fundamental forms} 

If we widen our viewpoint a bit, we can obtain a beautiful representation of possible geometric primitives for 3D projective space  (\cite{lle86}, p. 11). In this triangular arrangement, points, lines, and planes have equal rights; each occupies one corner of the triangle.  Consult Fig. \ref{fig:eltri}. 
 To begin with they exist alone and for themselves: point \emph{qua} point, line \emph{qua} line, plane \emph{qua} plane.  
 
 \subsubsection{Relationship through incidence}
 The other entries in the \emph{elemental triangle} arise by considering incidence relationships.
 Two  primitives  of different types come into a relationship via their \emph{incidence} properties. For example, given a point $P$ and a plane $p$, does the point lie on the plane; equivalently, does the plane pass through the point? If so, $P$ and $p$ are said to be \emph{incident}. One can ask, given a  point $P$, which planes are incident with $P$? Or, which points are incident with a plane $p$? 

The other entries in the triangle are obtained by taking an element from one corner and \quot{dropping it} into another corner of the triangle and observing the resulting incidence relationships that arise. For example, drop a plane $p$ (lower left corner) into the space of points (lower right). Imagine that any point incident with $p$ \quot{lights up}. Which figure will be thereby illuminated?  In this case, all the points of $p$ will light up. The resulting composite object is called a \emph{point field}. It is shown on the lower edge of the triangle to the left of the space of points. Similar combinations lead to the other 5 forms on the boundary of the triangle. 

\subsubsection{Line-based forms}
Notice that the lower left three (plane, plane pencil, plane bundle) occur in plane-based algebras, the lower right three (point, point range, point field) occur in point-based algebras.  The upper three entries do not occur in the geometric algebra framework discussed above. But since for this discussion we are treating the line  as an indivisible element in its own right, we obtain two further 2D composite forms, the \emph{line field} and \emph{line bundle}, consisting of all the lines incident with a plane, resp., with a point. They are dual to each other.  
Finally, in the center of the triangle is a primitive combining all three types, the \emph{line pencil}: the set of all lines incident with a point and with a plane that is incident with the point. It is self-dual.

\subsubsection{Dimensions of the forms}
The 10 entries of the triangle have different dimensions.  The three corners: point space, line space, and plane space, have dimensions 3, 4, 3 resp. The line and point fields, and the plane and line bundles, have dimension 2. The point range, line and plane pencil have dimension 1.  By duality, all 2-dimensional forms (fields and bundles) are structurally identical from the point of view of projective geometry, as are all 1-dimensional forms. They are sometimes referred to as as \emph{2D-}, resp., \emph{1D-}\emph{extents}.

\subsubsection{Concluding remarks}
The resulting arrangement is called the \emph{elemental triangle} and the entries are called the 10 \emph{fundamental forms} of $\RP{3}$, projective 3-space.  
In contrast to traditional geometry, this approach gives equal rights to each of the three fundamental entities point -- line -- plane. With these entities alone however no geometry is possible. That arises only when the seven composite forms are created out of incidence relationships among these three.  In this framework, what traditional geometry calls "point -- line -- plane" are just the three entries in the lower right corner of the triangle: point \emph{qua} point  -- point range (line of points) -- point field (plane of points) . Through the magic of geometric duality, this limited vocabulary is extended to the full decad depicted in the elemental triangle. A whole new world of geometry opens up with this extended vocabulary. This article can only touch on the possibilities. 


\subsection{The duck test for double algebra representations} 
The \emph{duck test} is used to show the equivalence of two entities based on their behaviors:
\begin{quote}
    If it looks like a duck, swims like a duck, and quacks like a duck, then it probably is a duck.
\end{quote}
In our case, we have shown that a double algebra approach supports geometric duality; we now argue that a system that supports geometric duality is in fact a double algebra representation.  Indeed, the foundation for the double algebra representation is the dual coordinate map $J: G \leftrightarrow G^*$ that maps a geometric entity to its dual representation. On the other hand, if I support geometric duality, then I can map a geometric entity $X$ to its dual form.  But that is equivalent to knowing the dual coordinate $J$ and hence I have a double algebra representation.  We conclude:

\begin{tcolorbox}[colback=ltltgreen]
\centering Support for geometric duality $\Leftrightarrow$ support for double algebra representation
\end{tcolorbox}

\subsection{Geometric duality and the real world}

Geometric duality can be seen and appreciated as a deep principle within the realm of geometry.  But it's not only a mathematical wonder.  The distinction between point-based and plane-based primitives appears also within natural science.  We turn now to indicate two sorts of such appearances. The first comes from the realm of classical euclidean rigid body mechanics. The second moves beyond euclidean space (our $G$ algebra) and brings in dual euclidean space ($G^*$ here) as a tool for scientific investigation.

\subsubsection{Example: rigid body kinematics and dynamics in euclidean PGA}
\label{sec:kmepga}
 It may come as a surprise to find out that geometric duality plays a key role in euclidean rigid body motion. 
 In 3D EPGA, both velocities and momenta, accelerations and forces, are all represented by bivectors (see \cite{gunn2011} for an introductory treatment). These kinematic and dynamical bivectors can be distinguished on the basis of geometric duality.
 
 \myboldhead{Kinematics is plane-based} The kinematic entities (velocities and accelerations) arise \textit{naturally} as bivectors in $G$, that is, as \textit{intersections of planes}. For example, a rotation is the geometric product of reflections in two planes that intersect in the axis of rotation and meet at half the angle of the rotation, and the logarithm of such a rotation is a scalar multiple of the bivector part.  As mentioned above, such a plane-wise line is called an \textit{axis}.  See Fig. \ref{fig:spearaxis}, left.

\myboldhead{Dynamics is point-based} 
Momenta and forces, on the other hand, are related to velocities and accelerations by the action of the inertia tensor. The inertia tensor $A: \bigwedge^2 \rightarrow (\bigwedge^2)^* $ is symmetric bilinear form that maps plane-based kinematic vectors (in $G$) that to point-based dynamic bivectors  in $G^*$.  Its inverse $A^{-1}$ does the reverse. A momentum or force is then \textit{naturally} a spear.  See Fig. \ref{fig:spearaxis}, right. You can also verify this in your own experience: if someone pushes on you with their finger, you can feel the spear, the line of force in its point-wise aspect, exerting directional pressure on your chest. This is the line in its spear aspect.

A further nice feature of the double algebra approach is that it is consistent with the traditional representation of the inertia tensor as a diagonal operator $G\rightarrow G^*$. In a single-algebra approach the distinction between spears and axes does not exist; every bivector is an axis and the matrix form is anti-diagonal. 

\myboldhead{Duality and outermorphisms} Technically, the inertia tensor is implemented within PGA as an \emph{outermorphism} $L: G\rightarrow G^*$. In the language of projective geometry, such a map is called a \emph{correlation} (since it maps points to planes), while a map $L: G\rightarrow G$ is called a \emph{collineation}.  Extending the standard concept of outermorphism, which is based on collineations, to include correlations is straightforward to achieve within the double algebra since the dual algebra $G^*$ is directly accessible via $\PN$. 

\myboldhead{Historical remarks} This relation of rigid body mechanics to 3D projective line geometry were established by Eduard Study in his ground-breaking work  \cite{study03}. He gave the two types of bivectors the suggestive names \textit{kineme} and \textit{dyname}, and also introduced the terms \emph{spear} and \emph{axis}.  

For a fuller discussion of rigid body mechanics in (euclidean and non-euclidean) PGA see  (\cite{gunnthesis}, Ch. 8 and 9.

\subsubsection{Dual euclidean space} \label{sec:des}

The dual algebra $G^*$ can be considered in its own right, as the point-based algebra $\pdclal{n}{0}{1}$. In this respect, it models the metric relationships in \emph{dual euclidean space} (DES).  An introduction to DES and scientific applications can be found in Ch. 10 of \cite{gunnthesis}. DES is obtained by dualizing euclidean space (ES) within projective space. Example: The absolute quadric of ES is the ideal plane and all its points and lines, DES features a single ideal point and all its planes and lines. If we think of ES as the result of expanding a sphere until it becomes a flat plane, DES arises in the opposite process, whereby the space contracts until all curvature is concentrated in a point. 

In his habilitation speech (1855) Bernard Riemannian announced his ground-breaking "empirical space" program (1855), whereby the metric properties of space are to be determined by empirical observation. DES fills a hole in the \quot{circle of geometries} available for this empirical determination. The classical geometries of constant curvature (hyperbolic, euclidean, elliptic) are traditionally located on a real line parametrized by their curvature. By adding DEG (corresponding to curvature $\infty$), this line can be bent to a closed, self-dual circle (the "one-point compactification" of the line) that opens up a whole new world of metric relationships available to scientific inquiry.

This chapter in \cite{gunnthesis} includes an overview of research in the last hundred years devoted to  preliminary results of the applications of DES in this way in the last 100 years. See for example \cite{adams80}, \cite{adams96}, \cite{lehrs2014} (especially Ch. 9), \cite{thomas09}, \cite{diener2021}. In almost all of this research, DES is considered as an active partner of ES, neither alone seems to provide the necessary framework.  This naturally lends itself to a double algebra approach where ES and DES are \emph{equal citizens} and represents an important infrastructure tool for this research agenda. 

In short, having direct access to the point-based as well as plane-based algebra serves not just mathematical and pedagogical goals, but also important scientific ones. 

\subsection{Geometric duality: Summing up}
Here are the chief results of the above considerations:
\begin{enumerate}
\item Dimension isn't inherent in the geometric primitive but depends on how space is conceptualized.
\item In a point-based algebra, points are 0-dimensional and planes are 2-dimensional (a point field); in a plane-based algebra, planes are 0-dimensional and points are 2-dimensional (a plane bundle).
\item The projective geometric principle of duality provides a logically rigorous and historically well-established basis for this dual conceptualization of space.
\item Every geometric primitive has two representations: a plane-based one and a point-based one.
\item \poincare duality $\PN$ maps between these two representations. We call the result a \emph{double-algebra representation}, Hodge duality $\CR$ results in a \emph{single algebra} approach. 
\item Support for geometric duality is equivalent to having a double algebra representation.  
\item Geometric duality is built into the structure of the physical world.
\end{enumerate}

\section{A bit better: Implementing the double algebra approach}
\label{sec:impdaap}
We come now to the end of our discussion. We began with the multiplicity of different \quot{duality maps} used in geometric algebra and ended up focusing our attention on two candidates that work for degenenerate metrics, \poincare and Hodge duality.  Sec. \ref{sec:compare} established that the coordinates obtained for the maps $\PN$ and $H$ are identical (and agree with $P$, in the case that the metric is non-degenerate). The only difference remaining is how the coordinates are interpreted. This led to the distinction between a double algebra representation (\poincare) versus a single algebra one (Hodge). We introduced the concept of geometric duality and established it is equivalent to the double algebra representation.

In this section we want to sketch a simple implementation of the double algebra approach. Once we have established it, we show how the implementation also has advantages for the single algebra approach. Hence the proposal can be seen as a single duality-neutral solution with advantages for the whole community.

This section is a placeholder for a treatment that takes into account the concrete details of modern software design. It's included here since the first step in planning a trip is to clarify the destination.


In this sense the duality maps we have examined share a common root. We denote it by $\mb{D}$ and leave open whether it represents $\PN$ or $\CR$ duality.  

The key requirement for supporting geometric duality is a mark showing whether a given multivector is in $G$ or in $G^*$. 
To keep track of this, add a boolean field \texttt{isDual} to the multi-vector instance with default value  \texttt{false}; in euclidean PGA, that means that all multi-vectors start out in $G$, the plane-based algebra. In dual euclidean PGA $\pclal{n}{0}{1}$, the default value would be \texttt{true}, since this is a point-based algebra. In their return value, $\mb{D}$ and $\mb{D}^{-1}$ flip the value of the \texttt{isDual} bit from the argument value, indicating that argument has been mapped to its geometrically dual form.


In fact we can use the \texttt{isDual} field to simplify the API for $\mb{D}$.   The idea is simple: $\mb{D}(x)$ is only valid when \texttt{x.isDual==false}, while $\mb{D}^{-1}$ only applies when \texttt{x.isDual==true}. Denoting the existing maps as $D_o$ and $DI_o$, resp., we can express the new single duality map $\mb{D}$ as follows:
\begin{lstlisting}
D(x) := {
    var dx = (!x.isDual) ? Do(x) : DIo(x); 
    dx.isDual = !x.isDual; 
    return dx;
}
\end{lstlisting}
This simplifies the API by replacing the two versions of the duality map (forward and back) with a single map that uses the \texttt{isDual} field to decide whether to invoke the $G\rightarrow G^*$ or the $G^*\rightarrow G$ version.

\subsection{Safe and correct duality}
We've seen that in the double algebra approach, calls to $\mb{D}$ and $\mb{D}^{-1}$ alternate.  
A glance at the three formulas for the regressive product in Sec. \ref{sec:compare} show that in these formulas the same is true for single algebra duality. 
 In all these formulas, $\mb{D}$ is invoked on both arguments, a product is applied, and $\mb{D}^{-1}$ is then applied to the result. Correct results are obtained only by \emph{alternating} calls to $\mb{D}$ and $\mb{D}^{-1}$. I am aware of no valid use case involving $\mb{D}$ and $\mb{D}^{-1}$ in which this principle of alternation is not respected. 
The \texttt{isDual} field provides a simple means to guarantee this alternation.

\myboldhead{Takeaway} Even users who see no direct need for a double algebra representation will profit from the \texttt{isDual} field since it simplifies the API (as seen above) and guarantees correct alternation in calls to $\mb{D}$ and $\mb{D}^{-1}$.



\subsection{Technical aspects of the \texttt{isDual} bitfield}

The introduction of a new field into a data structure has to be checked and integrated with the existing fields and functionality.  We have shown above how it is to be used and set when the duality map $\mb{D}$ is called. 

How does the \texttt{isDual} field interact with algebra products? We focus on the geometric product since others are derived from it.  Suppose $a$ and $b$ are two multi-vectors. Then we could write pseudo-code for the geometric product when the \texttt{isDual} fields are the same:
\begin{lstlisting}firstline=1,lastline=5]
gp(a,b) := {
    (a.isDual == b.isDual) ? 
    {var ab = a*b; ab.isDual = a.isDual; return ab;} :
    {??}
}
\end{lstlisting}

Here the \texttt{a*b} is the geometric product in the respective algebra (G when \texttt{isDual = false}, $G^*$ when  \texttt{isDual = true}. But when the \texttt{isDual} fields differ (indicated by the \texttt{??} text), what should happen?

\subsection{When the \texttt{isDual} fields don't match} 
\label{sec:hexc}
Exceptions are in the mind of the beholder. What one programmer considers an error may be for another a feature. This leads to two strategies for dealing with exceptional values of the \texttt{isDual} field.  

\textbf{Strict interpretation:} When the geometric product is called on two multi-vectors whose \texttt{isDual} fields differ, the geometric product is not defined. In the double algebra representation, this is self-evident. And in the single algebra interpretation, we  pointed out above that the proper alternation of the \texttt{isDual} field is important to obtain correct answers. 
 Conclusion: this exception should signal an error. Such a strategy can be  useful for programmers to improve their code.

\textbf{Lenient interpretation:} It's not an error, it's a type conversion challenge. Type conversion is well-known and understood topic in numerical mathematics.  Consider the product of two real numbers $a*b$.  Each of $a$ and $b$ has a specific type: \texttt{int, float, double,...}.  At the hardware level, multiplication is only defined for arguments of the same type. In order to handle all possible combinations of types, a hierarchy of promotions has been developed and embedded in compilers so that the output value has optimal properties, in this case, retains the maximal precision. If for example \texttt{int a = 5; double b = 10; a*b;} will implicitly convert $a$ to a \texttt{double} value before calculating $a*b$, which then is of type \texttt{double}. 

The double algebra approach invites a similar treatment.   $\PN$ maps between two different representations of the same geometric element, one plane-based, one point-base, analogous to having float and double representations of a real number. In this analogy, when $a$ and $b$ are in different algebras, we establish a type promotion rule to obtain elements of the same algebra before carrying out the operation. When we are doing euclidean geometry, the primary algebra of interest is $G$, so we promote the element of $G^*$ to $G$ before carrying out the calculation. The complete pseudo-code becomes:
\begin{lstlisting}
gp(a,b) := {
    if (a.isDual != b.isDual)  {a.isDual ? a = D(a) : b = D(b);}  
    var ab = a*b; ab.isDual = a.isDual; return ab; 
}
\end{lstlisting}
Note that this treatment only applies when the \texttt{isDual} fields of the original two arguments don't match. If the arguments match, the multiplication takes place in the algebra implied by the \texttt{isDual} fields. 

Similar remarks apply to handling exceptions raised when calling the unified duality map $\mb{D}$ introduced above.  In all cases there is a strict and a lenient interpretation.

In general, users should be able to control how the \texttt{isDual} field is configured and how associated exceptions are handled.

\subsection{Summing up}
To sum up, introducing the \texttt{isDual} field to the multi-vector class yields the following valuable features:
\begin{enumerate}
    \item It provides checks that the proper alternation of calls to $\mb{D}$ and $\mb{D}^{-1}$ are observed.
    \item It allows unifying the duality map and its inverse into a single map.
    \item It identifies whether a multi-vector is in $G$ or $G^*$.  For example, visualization of a simple bivector can take into account the difference mentioned above between \emph{spears} (lines in a point-based algebra) and \emph{axes} (lines in a plane-based algebra). 
    \item This also allows the fundamental duality noted above between kinematics ($G$) and dynamics ($G^*$) to be directly mirrored in the underlying computational framework.
\end{enumerate}

\section{Conclusions}
\label{sec:conc}
In geometric algebras with degenerate metrics, such as the algebra $\pdclal{n}{0}{1}$ for euclidean space $\mb{E}^{n}$, non-metric methods are required in order to carry our important tasks that traditionally have been carried out using a non-degenerate pseudoscalar, such as coordinate computation and regressive product.

We have focused on two grade-reversing maps to replace multiplication by the pseudoscalar: the dual coordinate map $\PN: G \rightarrow G^*$ and the right complement map $\CR: G \rightarrow G$. In fact, we showed that in terms of a standard basis, the two maps take identical forms. The only difference between $\PN$ and $\CR$ is one of interpretation: whether the resulting coordinate index set should be interpreted as superscripts (in $G^*$) or as subscripts (in $G$).  

This isomorphism hides an important difference. Only $\PN$ is coordinate-free. This follows from the decomposition of $\PN = \RL\circ \CR$ into right complement ($\CR$) followed by left complement ($\RL$), yielding the well-known fact that $\PN$ is an identity map on geometric primitives, and hence a \emph{natural} map.   This also is consistent with the traditional name, \emph{dual coordinate map}, for this form of duality.

 We introduced the concept of \emph{geometric duality} to capture the double nature of all geometric primitives, appearing once in the plane-based and once in the point-based algebra. We discussed the scientific examples, such as rigid body mechanics, where geometric duality plays a crucial role. We saw that geometric duality naturally maps onto a \emph{double algebra} representation (e. g., \poincare duality) but not onto a \emph{single algebra} representation (e.g., Hodge duality). 

Based on these observations, we proposed a software design designed in the first place to produce a double algebra representation, but which also has advantages for the single algebra approach. It introduces a bit field in the multi-vector data structure marking whether that multi-vector is an element of $G$ or $G^*$. Calls to the duality function $\mb{D}$ flip that bit. We also showed that exceptions related to such a bit field can be handled using promotion rules modeled on the familiar promotion rules used for multiplying numbers of different data types together, rather than treating such exceptions as errors.


The goal is creating software frameworks for PGA that are \quot{a bit better}, providing support for all members of the  community regardless of their duality affiliation. This would produce a welcome \quot{unity in duality}.


\bibliography{GunnRef}
\bibliographystyle{alpha}

\end{document}